\documentclass{amsart}
\usepackage{graphicx}
\usepackage{amsmath}
\usepackage{amscd}
\usepackage{amsfonts}
\usepackage{amssymb}

\numberwithin{equation}{section}

\newcommand{\be}{\begin{equation}}
\newcommand{\ee}{\end{equation}}
\newcommand{\bes}{\begin{equation*}}
\newcommand{\ees}{\end{equation*}}

\newcommand{\mb}[1]{\mathbb{#1}}

\begin{document}

\title[Three remarks]{Three remarks on a question of Acz\'{e}l}

\author{Orr Moshe Shalit}
\address{Department of Pure Mathematics, University of Waterloo.}
\email{oshalit@math.uwaterloo.ca}

\maketitle

Consider the functional equation
\be\label{eq:fe}
f(x^2 R) = \frac{k}{2xR}f(x) \,\, , \,\, x > \frac{1}{R} ,
\ee
where $R,k$ are positive constants.
At the 49th International Symposium on Functional Equations 2009, J. Acz\'{e}l \cite[p. 195]{report} presented the general solution of this equation as
\be\label{eq:gs}
f(x) = \frac{(\ln(xR))^c}{xR}p(\log_2(\ln(xR))) , 
\ee
where $c = \log_2(k/2)$ and $p$ is an arbitrary periodic function of period $1$  on $\mb{R}$ (for proof define $p(s) = (\frac{2}{k})^s \exp(2^s) f\left(\frac{\exp(s^2)}{R}\right)$), and then asked for the \emph{monotonic} solutions.


A special solution for the functional equation (\ref{eq:fe}) is
\be\label{eq:special}
\varphi_c(x) = \frac{(\ln(xR))^c}{xR},
\ee
obtained by taking $p \equiv 1$ in (\ref{eq:gs}). This function is monotonic precisely when $c \leq 0$, that is, when $k \leq 2$. 
Are all monotonic solutions of (\ref{eq:fe}) scalar multiples of $\varphi_c$? We show that when $k=2$ the answer is \emph{yes}, and that it is \emph{no} when $k<2$. We then consider the effects of imposing continuity and differentiability conditions at $1/R$. Here, too, the results depend on the value of $k$.

\section{Monotonicity}

First, assume that $k = 2$. We claim that a function of the form 
\bes
\frac{1}{xR}\cdot  p(\log_2(\ln(xR)))
\ees
is monotonic only if $p \equiv const$. Indeed, if $p$ is not constant, then it takes two distinct values $M>m$. For any $x_1, x_2$ sufficiently close to $1/R$, 
\bes
M \frac{1}{x_1 R} > m \frac{1}{x_2 R} .
\ees
Let $x_1$ be sufficiently close to $1/R$ which satisfies $p(\log_2(\ln(x_1 R))) = M$, and let $x_2 > 1/R$ be smaller than $x_1$ such that 
$p(\log_2(\ln(x_2 R))) = m$. Then $x_2 < x_1$, but $f(x_1) > f(x_2)$, so $f$ is not monotonic decreasing. If $f$ is an increasing function then $-f$ is a decreasing function, so we conclude that the only monotonic solutions to (\ref{eq:fe}) are $f(x) = \frac{\lambda}{xR} = \lambda \varphi_0(x)$, $\lambda \in \mb{R}$.  

Now we consider the case $k<2$. We claim that there are infinitely many differentiable functions $p$ such that (\ref{eq:gs}) is a monotonic solution. Indeed, differentiating (\ref{eq:gs}) we find
\be\label{eq:f'}
f'(x) =\frac{\ln(xR)^{c-1}}{x^2 R}\Big((c-\ln(xR))\cdot p(\log_2(\ln(xR))) + \log_2  e \cdot  p'(\log_2(\ln(xR))) \Big).
\ee
Keeping in mind that $c<0$ in this case, it is easy to see that, so long as $p$ is bounded away from $0$ and $|p'|$ is bounded from above by a small enough number, $f'$ has a constant sign in $(1/R,\infty)$. Thus, there are many monotonic - even differentiable - solutions other than $\varphi_c$.

\section{Continuity at $1/R$}


If the domain in (\ref{eq:fe}) is changed to $x \in (0,1/R)$, then the general solutions is
\be\label{eq:gs2}
f(x) = \frac{(-\ln(xR))^c}{xR}p(\log_2(-\ln(xR))) .
\ee

Note that $f$ can be extended to a continuous function on $[1/R,\infty)$ that satisfies (\ref{eq:fe}) at $x = 1/R$ only if $k=2$ or $f(1/R) = 0$. None of the solutions for $k<2$ can be extended to a continuous solution on $[1/R,\infty)$. 

On the other hand, every solution for $k>2$ can be extended to $x=1/R$ by defining $f(1/R)=0$. Pasting these solutions with (\ref{eq:gs2}) for $0<x<1/R$, we see that there are many continuous solutions for (\ref{eq:fe}) when $k>2$ on the interval $(0,\infty)$. All of these solutions vanish at $0$.
 
Finally, when $k=2$, we see that for every $\lambda \in \mb{R}$, the function $\lambda \varphi_0$ is the unique continuous solution for (\ref{eq:fe}) that satisfies $f(1/R) = \lambda$. 
Of course, this is a solution of (\ref{eq:fe}) on the entire interval $(0,\infty)$.

\section{Continuous differentiability}

The only case worth discussing in the setting of contionuous differentiablity is $k>2$. In this case, we saw that the functional equation has continuous solutions on the the interval $(0,\infty)$ given by (\ref{eq:gs}) and (\ref{eq:gs2}). We ask: which of these solutions is continuously differentiable? Examining (\ref{eq:gs}), we see that we may limit the discussion to solutions of the form (\ref{eq:gs}) with $p$ continuously differentiable. 

When $k \in (2,4]$, so $c\in (0,1]$, we see by examining (\ref{eq:f'}) that $f'$ will have a limit at $x = 1/R$ only if 
\bes
\lim_{x \rightarrow 1/R^+} (c-\ln(xR))\cdot p(\log_2(\ln(xR))) + \log_2  e \cdot  p'(\log_2(\ln(xR))) = L ,
\ees
and $L$ must $0$ when $c<1$. But since $p$ is periodic, this is possible only if $p$ satisfies the differential equation
\bes
p' + \frac{c}{\log_2 e} p = {L}/{\log_2 e} .
\ees
But the only periodic solutions of this differential equation are constants, thus any continuously differentialbe solution when $k \in (2,4]$ is multiple of $\varphi_c$.

When $k>4$ any continuously differentiable periodic $p$ will give rise to a continuously differentiable solution.


\bibliographystyle{amsplain}

\end{document}